\documentclass[11pt]{article}
\usepackage{amssymb}

\title{On the consistency strength of the Inner Model Hypothesis
} 
\author{Sy-David Friedman (KGRC, Vienna)\footnote{The first author
     was supported by FWF Grants P16334-NO5 and P16790-NO4.}\\Philip Welch
   (Bristol)\\Hugh
   Woodin (Berkeley)
}
\date{10 August 2006}
\newcommand{\nod}{\noindent}

\newtheorem{thm}{Theorem}
\newtheorem{lem}[thm]{Lemma}

\newtheorem{prop}[thm]{Proposition}
\newcommand{\lan}{\langle}
\newcommand{\ran}{\rangle}

\newcommand{\card}{\hbox{card }}
\newcommand{\hcard}{\hbox{hcard }}

\newcommand{\Ord}{\hbox{Ord}}

\newcommand{\bs}{\bigskip}
\newcommand{\noi}{\noindent}
\newcommand{\cof}{\hbox{cof }}
\newcommand{\Cof}{\hbox{Cof }}

\newcommand{\Lim}{\hbox{Lim }}
\newcommand{\ot}{\hbox{ot }}

\newcommand{\ordertype}{\hbox{ordertype }}

\begin{document}
\maketitle

The \emph{Inner Model Hypothesis (IMH)} and the \emph{Strong Inner Model
Hypothesis (SIMH)} were introduced in \cite{imh}. In this article we
establish some upper and lower bounds for their consistency strength.

\bs

We repeat the statement of the IMH, as presented in \cite{imh}.
A sentence in the language of set theory is \emph{internally
  consistent} iff it holds in some 
(not necessarily proper) inner model. The meaning of internal
consistency depends on 
what inner models exist: If we enlarge the universe, it is
possible that more statements become internally consistent. 
The \emph{Inner Model Hypothesis} asserts that the 
universe has 
been maximised with respect to internal consistency:

\bs

\noi
\emph{The Inner Model Hypothesis (IMH):} 
If a statement $\varphi$ without parameters holds in an inner model of
some
outer model of $V$ (i.e., in some model compatible with $V$), then it
already holds in some inner model of $V$.

\bs

\noi
Equivalently: If $\varphi$
is internally consistent in some outer model of $V$ then it is already
internally consistent in $V$.
This is formalised as follows. Regard $V$ as a
countable model of G\"odel-Bernays class theory,
endowed with countably many sets and classes. Suppose that $V^*$
is another such model, with the same ordinals as $V$. Then $V^*$ is an
\emph{outer model  of $V$} ($V$ is an \emph{inner model of $V^*$}) iff
the sets of $V^*$ include the
sets of $V$ and the classes of $V^*$ include the classes of $V$.
$V^*$ is \emph{compatible with $V$} iff $V$ and $V^*$ have a common
outer model.

\bs
\noi
\emph{Remark.} 
The Inner Model Hypothesis, like L\'evy-Shoenfield
absoluteness, is a form of absoluteness 
between $V$ and arbitrary outer models of $V$, which need not
be generic extensions of $V$. Formally speaking,
the notion of ``arbitrary outer model'' does depend on the background
universe in which $V$ is situated as a countable model. However, 
a typical model of the IMH is minimal in the sense
that for some real $R$, it is the smallest transitive model of
G\"odel-Bernays containing $R$ (see Theorem \ref{woodin2} below). 
For minimal models, the choice of background universe is irrelevant,
and if there is a model of the IMH then there is a minimal one. Thus
we may in fact regard the IMH as an intrinsic hypothesis about $V$,
independent of any background universe. 
An alternative way to
``internalise'' the IMH is to restrict the notion of outer model to
class-generic extensions which preserve the axioms of
G\"odel-Bernays. This weakened form of the IMH is expressible in the
language of class theory, and the results of this paper would not be
affected by this change. However, this is not in the spirit of
\cite{imh}, where the IMH is introduced as a form of
absoluteness which is no way limited by the notion of forcing.

\begin{thm} \label{imh} (\cite{imh})
The Inner Model Hypothesis implies that for some real $R$,
ZFC fails in $L_\alpha[R]$ for all ordinals $\alpha$. In particular,
there are no inaccessible cardinals, the reals are not closed under
$\#$ and the singular cardinal hypothesis holds.
\end{thm}

\begin{thm} \label{meas}
The IMH implies that there is an inner model with measurable cardinals
of arbitrarily large Mitchell order.
\end{thm}

\noi
\emph{Proof.} Assume not and let $K$ denote Mitchell's core model for
sequences of 
measures (see \cite{mitchell}). Let $\delta$ be the maximum of
$\omega_1$ and the strict supremum of the Mitchell orders of measurable
cardinals in $K$. By Mitchell's Covering Theorem for
$K$ we have:

\bs

\noi
\emph{$(*)$ $\cof(\alpha)\ge\delta$, $\alpha$ regular in $K\rightarrow
\cof(\alpha)=\card(\alpha)$.}

\bs

\noi
Now iterate $K$ by applying each measure of order $0$ exactly
once, i.e., if $K_i$ is the $i$-th iterate of $K$, $K_{i+1}$ is formed
by applying the measure of order $0$ in $K_i$ at $\kappa_i$, where 
$\kappa_i$ is the least measurable of $K_i$ not of the form
$\pi_{i_0i}(\kappa_{i_0})$ for any $i_0<i$. It is easy to see that
$i<j$ implies $\kappa_i<\kappa_j$ and hence this iteration is normal.
Let $\sigma_{ij}:K_i\to K_j$ be the 
resulting iteration map from $K_i$ to $K_j$ for $i\le j\le\infty$, 
let $K'$ denote $K_\infty$ and
let $\sigma:K\to K'$ denote $\sigma_{0\infty}$. 

\begin{lem} \label{k-prime-analysis}
(a) For any $\alpha$ and any $i$, there is at most one $j\ge i$ such
that 
$\sigma_{j,j+1}$ is discontinuous at $\sigma_{ij}(\alpha)$. If there
is no such $j\ge i$, then $\sigma_{i\infty}$ is continuous at
$\alpha$.\\
(b) For any $\alpha$ and any $i$, the $K_i$-cardinality of
$\sigma_{i\infty}(\alpha)$ is at most $\alpha^+$ of $K_i$.
If $\sigma_{i\infty}$ is continuous at $\alpha$ and $\alpha$ is a
$K_i$-cardinal, then the $K_i$-cardinality of 
$\sigma_{i\infty}(\alpha)$ equals $\alpha$.\\
(c) If $\alpha$ is a limit cardinal of $K_i$ then $\alpha$ is a
closure point of $\sigma_{i\infty}$, i.e.,
$\sigma_{i\infty}[\alpha]\subseteq\alpha$. If 
in addition the $K_i$-cofinality of $\alpha$ is not measurable in
$K_i$, then $\alpha$ is a fixed point of $\sigma_{i\infty}$.
\end{lem}

\noi
\emph{Proof.} (a) The map $\sigma_{j,j+1}$ is discontinuous at 
$\beta$ iff $\beta$ has $K_j$-cofinality $\kappa_j$. Thus if
the $K_j$-cofinality of $\sigma_{ij}(\alpha)$ is not $\kappa_j$ for any
$j\ge i$, it follows that $\sigma_{j,j+1}$ is continuous at
$\sigma_{ij}(\alpha)$ for all $j\ge i$. Otherwise let $j$ be least so
that 
$\sigma_{ij}(\alpha)$ has 
$K_j$-cofinality $\kappa_j$; then for $k$ greater than $j$,
$\sigma_{ik}(\alpha)$ has $K_k$-cofinality $\sigma_{jk}(\kappa_j)$,
which 
by definition of the $\kappa_k$'s does not equal $\kappa_k$. It
follows that $\sigma_{k,k+1}$ is continuous at $\sigma_{ik}(\alpha)$
for $k$ greater than $j$, as desired. The last statement is
immediate.\\ 
(b) Define the ordering $\prec$ as follows: $(j,\alpha)\prec
(k,\beta)$ iff $j\ge k$ and $\alpha<\sigma_{kj}(\beta)$. The
relation $\prec$ is a well-founded partial ordering. We prove the
desired property of $(i,\alpha)$ by induction on $\prec$. 

We may assume that $\alpha$ is a cardinal of $K_i$. 
If $\sigma_{ij}(\alpha)$ does not have $K_j$-cofinality $\kappa_j$ for
any $j\ge i$, then $\sigma_{i\infty}$ is continuous at $\alpha$ and
by induction applied to pairs $(i,\bar\alpha)$, $\bar\alpha<\alpha$,
we 
have that 
the $K_i$-cardinality of $\sigma_{i\infty}(\alpha)=\sup
\sigma_{i\infty}[\alpha]$ is equal to 
$\alpha$. Otherwise, choose the unique $j\ge i$ so that
$\sigma_{ij}(\alpha)$ has $K_j$-cofinality $\kappa_j$; then
$\sigma_{ij}$ is continuous at $\alpha$ and therefore by induction
applied to pairs $(i,\bar\alpha)$, $\bar\alpha<\alpha$, 
$\sigma_{ij}(\alpha)$ has $K_i$-cardinality $\alpha$. 
Let $\alpha^*$ denote $\sigma_{ij}(\alpha)$.
Now $\sigma_{i,j+1}(\alpha)=\sigma_{j,j+1}(\alpha^*)$ has
$K_j$-cardinality $(\alpha^*)^+$ of $K_j$. And by induction applied 
to pairs $(j+1,\beta)$, $\beta< \sigma_{j,j+1}(\alpha^*)$, we have
that  $\sigma_{j+1,\infty}(\beta)$ has $K_{j+1}$-cardinality at most
$\beta^+$ of 
$K_{j+1}$ for each $\beta<\sigma_{j,j+1}(\alpha^*)$.
So as $\sigma_{j+1,\infty}$ is continuous at
$\sigma_{j,j+1}(\alpha^*)$ and $\sigma_{j,j+1}(\alpha^*)$ is a
cardinal of $K_{j+1}$, it follows 
that the $K_{j+1}$-cardinality of $\sigma_{j+1,\infty}
(\sigma_{j,j+1}(\alpha^*))=\sigma_{j\infty}(\alpha^*)$ is
$\sigma_{j,j+1}(\alpha^*)$, and therefore the $K_j$-cardinality 
of $\sigma_{j\infty}(\alpha^*)$ is 
$(\alpha^*)^+$ of $K_j$. 
As $\alpha^*$ has $K_i$-cardinality $\alpha$, it follows that 
the $K_i$-cardinality of
$\sigma_{i\infty}(\alpha)=\sigma_{j\infty}(\alpha^*)$ is at most 
$\alpha^+$ of $K_i$, as desired.\\
(c) The first statement follows immediately from (b). For the second
statement, suppose that $\alpha$ is a closure point of
$\sigma_{i\infty}$ and the $K_i$-cofinality of $\alpha$ is not
measurable in $K_i$. We show by induction on $j\ge i$ that 
$\sigma_{ij}$ is continuous at $\alpha$, and therefore that $\alpha$
is a fixed point of $\sigma_{ij}$: This is vacuous if $j=i$. If
$\alpha$ is a fixed point of $\sigma_{ij}$ then the $K_j$-cofinality
of $\alpha$ is not measurable in $K_j$ by elementarity, and therefore
$\sigma_{j,j+1}$ is continuous at $\alpha$; it follows that
$\sigma_{i,j+1}$ is also continuous at $\alpha$. For limit $j$, the
continuity of $\sigma_{ij}$ at $\alpha$ follows from the continuity of
the $\sigma_{ik}$ at $\alpha$ for $k<j$.
$\Box$

\begin{lem} \label{k-prime-covering}
$(*)$ holds with $K$ replaced by $K'$.
\end{lem}

\noi
\emph{Proof.} It suffices to show by induction on $i$ that $(*)_i$
holds, where $(*)_i$ is $(*)$ with $K$ replaced by $K_i$.

\bs

\noi
Base case: $(*)_0$ is just $(*)$.

\bs

\noi
Successor case: 
Suppose that $(*)_i$ holds and that $\alpha$ is
$K_{i+1}$-regular with 
cofinality at least $\delta$.
We may assume that $\alpha$ is greater than $\kappa_i$,
else $\alpha$ is $K_i$-regular and we are done by induction. 
If $\alpha$ is at most $\sigma_{i,i+1}(\kappa_i)$ then $\alpha$ has 
$K_i$-cardinality $\kappa_i^+$ of $K_i$, and, as $K_{i+1}$ and
$K_i$ 
contain the same $\kappa_i$-sequences of ordinals, $\alpha$ has
$K_i$-cofinality $\kappa_i^+$ of $K_i$. So $\alpha$ and $\kappa_i^+$
of $K_i$ have the same cardinality and cofinality, so we are
done by induction.

\bs

Now suppose that $\alpha$ is greater than $\sigma_{i,i+1}(\kappa_i)$. 
Represent $\alpha$ in $K_{i+1}$, the ultrapower of $K_i$, by $[f]$
where 
$f:\kappa_i\to\Ord$. We may assume that $f$ is either constant or
increasing, and also that $f(\gamma)$ is $K_i$-regular and greater than
$\kappa_i$ for all $\gamma<\kappa_i$. 
If $f$ is constant then $\alpha=\sigma_{i,i+1}(\bar\alpha)$ for some
$\bar\alpha$ which is regular in $K_i$ and greater than $\kappa_i$;
then $\alpha$ and $\bar\alpha$ have the same cofinality and
cardinality, so we are done by induction. 
Thus we may assume that $f$ is increasing.

\bs

Now the 
$K_i$-cofinality of $\alpha$ is at least the supremum $\mu$ of the
$f(\gamma)$'s, as using the regularity of the $f(\gamma)$'s, we can
everywhere-dominate any set in $K_i$ of 
$f(\gamma)$-many 
functions from $\kappa_i$ into $\prod_{\gamma'>\gamma} f(\gamma')$ by a
single such function in $K_i$. As $\mu$ is $K_i$-singular, the
$K_i$-cofinality of $\alpha$ is in fact greater than $\mu$. 
And the $K_i$-cardinality of $\alpha$ is $\mu^{\kappa_i}=\mu^+$ of
$K_i$. It follows that $\alpha$ and $\mu^+$ of $K_i$ have the same
cofinality and cardinality, so we are done by induction.

\bs

\noi
Limit case:
Suppose that $i$ is a limit and $\alpha$ is $K_i$-regular with
cofinality at least $\delta$. By Lemma \ref{k-prime-analysis} (a), we
can choose $i_0<i$ such that 
$\alpha$ equals as $\sigma_{i_0i}(\bar\alpha)$, where $\bar\alpha$ is
regular in $K_{i_0}$ and $\sigma_{i_0i}$ is continuous at
$\bar\alpha$. It follows by Lemma \ref{k-prime-analysis} (a) that
$\alpha$ and $\bar\alpha$ have the same cardinality and cofinality,
and therefore we are done by induction. $\Box$

\begin{lem} \label{k-prime-sing}
If $\lambda$ is a limit cardinal then $\cof^{K'}(\lambda)$ is not
measurable in $K'$.
\end{lem}

\noi
\emph{Proof.} Let $\kappa$ denote the $K$-cofinality of $\lambda$.
If $\kappa$ is not measurable in $K$ then 
by Lemma \ref{k-prime-analysis} (c), 
$\lambda$ is a fixed point
of the iteration $\sigma:K\to K'$ and therefore the result follows by
elementarity. Otherwise we claim that $\kappa$ must equal $\kappa_i$
for some $i$: If not, then as $\kappa$ is a closure point of $\sigma$,
$\kappa$ would also be a fixed point of $\sigma$ and therefore
$\kappa$ is measurable in each $K_i$. By the definition of the
$\kappa_i$'s, it must be that for each $i$, $\kappa$ is either of the
form $\sigma_{i_0i}(\kappa_{i_0})$ for some $i_0<i$, or $\kappa_i$ is
less than $\kappa$. But for sufficiently large $i$, $\kappa_i$ cannot
be less than $\kappa$ and so $\kappa$ is of the form 
$\sigma_{i_0i}(\kappa_{i_0})$ for some $i_0<i$, which implies that
$\kappa$ equals $\kappa_{i_0}$, as $\kappa$ is a fixed point of 
$\sigma_{i_0i}$.

So choose $i$ so that $\kappa=\kappa_i$.
Then $\kappa$ and $\lambda$ are fixed points of $\sigma_{0i}$ 
and therefore $\lambda$ has $K_i$-cofinality $\kappa_i$.
As $K_{i+1}$ has the same $\kappa$-sequences as $K_i$, 
it follows that $\lambda$
has cofinality $\kappa$ in $K_{i+1}$ and therefore also in $K'$. But
since we applied the order $0$ measure on $\kappa$ to form $K_{i+1}$,
$\kappa$ is not measurable in $K_{i+1}$ and therefore not measurable
in $K'$, as desired. $\Box$  

\bs

Now we apply the technique of ``dropping along a square sequence''
from \cite{facts}. Define a function $d: \Ord\to\omega$ as 
follows. Fix a lightface $K'$-definable global $\Box$-sequence $\lan
C_\alpha\mid \alpha$ singular in $K'\ran$: $C_\alpha$ is closed
unbounded in $\alpha$ with ordertype less than $\alpha$ for each
$K'$-singular $\alpha$ and $C_{\bar\alpha}=C_\alpha\cap\bar\alpha$
whenever  $\bar\alpha$ is a limit point of $C_\alpha$. If $\alpha$ is
not $K'$-singular then $d(\alpha)=0$. Otherwise define:

\bs

\noi 
$\alpha_0=\alpha$\\
$\alpha_1=\ot(C_{\alpha_0})$\\
$\alpha_2=\ot(C_{\alpha_1})$\\
$\cdots$\\
$\alpha_{n+1}=\ot(C_{\alpha_n})$,

\bs

\noi
as long as $\alpha_n$ is $K'$-singular, and let $d(\alpha)$ be the
least $n$ such that $\alpha_n$ is not
$K'$-singular. $\alpha_{d(\alpha)}$ is the $K'$-cofinality of
$\alpha$.

\begin{lem} \label{main} (Main Lemma, after \cite{facts}) 
For each $n$ there is a ZFC-preserving class forcing $P_n$ that adds a
CUB class $C_n$ of singular cardinals such that for all $\alpha$ in
$C_n$ of cofinality at least $\delta$, $d(\alpha)$ is at least $n$.
\end{lem}

\noi
\emph{Proof.} We use the following.

\begin{lem} \label{drop}
Suppose $k<m$, $\alpha\ge\delta$, $\alpha$ is
regular and $C$ is a closed set of ordertype $\alpha^{+m} +1$, consisting of
ordinals $\ge$
$\alpha^{+m}$ (where $\alpha^{+0}=\alpha$,
$\alpha^{+(p+1)}=(\alpha^{+p})^+)$. Then $(C\cap\{\beta\mid
d(\beta)\ge k+1\})\cup\Cof(<\delta)$ has a closed subset of ordertype
$\alpha^{+(m-k-1)}+1$. 
\end{lem}

\noi
\emph{Proof.} The proof is by induction on $k$, using Lemma
\ref{k-prime-covering}.

\bs

\noi
Suppose $k=0$. Let $\beta$ be the
$\alpha^{+(m-1)}$-st element of $C$. Then $\beta$ is $K'$-singular
since its cofinality
($=\alpha^{+(m-1)}$) is at least $\delta$ and less than its
cardinality ($\ge\alpha^{+m}$). Similarly, each 
element of $\Lim (C\cap \beta)$ of cofinality $\ge\delta$ is
$K'$-singular and therefore 
$\Lim (C\cap \beta)$ is a closed subset of $(C\cap\{\beta\mid
d(\beta)\ge 1\})\cup\Cof(<\delta)$ of ordertype $\alpha^{+(m-1)}+1$,
as desired.

\bs

Suppose that the Lemma holds for $k$ and let $m+1>k+1$, $C$ a closed
set of ordertype $\alpha^{+(m+1)}+1$ consisting of ordinals
$\ge\alpha^{+(m+1)}$. Then $\mu=\max C$ is $K'$-singular, as its
cofinality is at least $\delta$ and less than its cardinality. 
Let $\beta$ be the 
$(\alpha^{+m}+\alpha^{+m})$-th element of $C\cap C_\mu$.
$\beta$ is $K'$-singular as its
cofinality is at least $\delta$ and less than its cardinality.
Let $\bar\beta$ be the 
$\alpha^{+m}$-th element of $C$. 
Then
$\bar C=\{\ot C_\gamma\mid \gamma\in C\cap\Lim
C_\beta\cap[\bar\beta,\beta]\}$ is a 
closed set of ordertype $\alpha^{+m}+1$ consisting of ordinals $\ge$
$\alpha^{+m}$. By induction there is a closed $\bar D$ contained
in $(\bar C\cap\{\gamma\mid d(\gamma)\ge k+1\})\cup\Cof(<\delta)$ of ordertype
$\alpha^{+(m-k-1)}+1$. But then $D=\{\gamma\in C\cap\Lim C_\beta\mid
\ot C_\gamma\in 
\bar D\}$ is a closed subset of $(C\cap\{\gamma\mid d(\gamma)\ge
k+2\})\cup\Cof(<\delta)$ 
of ordertype $\alpha^{+(m-k-1)}+1$. As $m-k-1=(m+1)-(k+1)-1$, we are
done. $\Box$ (Lemma \ref{drop})

\bs

Lemma \ref{main} now follows: Let $P_n$ consist of closed sets $c$ of
singular cardinals such that

\bs

\noi
$\alpha\in c$, $\cof(\alpha)\ge\delta\to d(\alpha)\ge n$,

\bs

\noi
ordered by end-extension. Lemma \ref{drop} implies that this forcing
is $\kappa$-distributive for every cardinal $\kappa$. $\Box$

\bs

Now for each $n$ there is an outer model of $V$ containing a real $R_n$ such
that in $L[R_n]$:

\bs

\noi
$(*)_{R_n}$ $R_n$ codes a CUB class $C_{R_n}$ of singular cardinals and an
iterable, universal extender model $K'_{R_n}$ such that\\
a. $d_{R_n}(\alpha)\ge n$ for $\alpha$ in $C_{R_n}$ of sufficiently large
cofinality, where $d_{R_n}(\alpha)$ is defined in $K'_{R_n}$ just like
$d(\alpha)$ is defined in $K'$.\\
b. $\alpha\in C_{R_n}\to \cof(\alpha)$ in $K'_{R_n}$ is not measurable in
$K'_{R_n}$.

\bs

\noi
This is because we can use Lemma \ref{main} to force a CUB class $C_n$
of singular cardinals such that $d(\alpha)\ge n$ for all $\alpha$ in $C_n$ of
sufficiently large cofinality, and then $L$-code the model $\lan
V, C_n, K'\ran$ by a real $R_n$. The extender model $K'$ is universal
in the extension as successors of strong limit cardinals are not
collapsed and therefore weak covering holds relative to $K'$ in the
extension at all such cardinals of sufficiently large cofinality.

\bs

Applying the IMH, there are such reals $R_n$ in $V$. As each $R_n$
codes a CUB class of singular cardinals, the $K$ of $L[R_n]$ is
universal and therefore so is the $K_{R_n}$ arising from
$(*)_{R_n}$. Now co-iterate the $K_{R_n}$'s to a single $K^*$,
resulting in embeddings $\pi_n:K_{R_n}\to K^*$. As
singular cardinals in $C_{R_n}$ of sufficiently large cofinality are
fixed by $\pi_n$ (as their $K_{R_n}$-cofinality is not measurable in
$K_{R_n}$), it follows that there is a single
$\gamma$ belonging to all of the $C_{R_n}$'s which is fixed by all of
the $\pi_n$'s. But then $d^*(\gamma)\ge n$ for each $n$, where
$d^*(\gamma)$ is defined relative to $K^*$ just like $d(\gamma)$ was
defined relative to $K'$. This is a contradiction. $\Box$

\bs

   For each real $x$ let $M_x$, if it exists,  be the minimum transitive
set  model of $ZFC$
containing $x$. Thus $M_x$ has the form $L_\mu[x]$ for some 
countable ordinal $\mu = \mu(x)$. If $d$ is a Turing degree we write
$M_d$, $\mu(d)$ for $M_x$, $\mu(x)$ ($x$ in $d$). 

\begin{thm} \label{woodin2}
Assume the existence of a Woodin cardinal with an
inaccessible above. Then the IMH is consistent. 
Moreover for all $d$ in a cone of Turing degrees, $M_d$ exists and
satisfies the IMH.
\end{thm}

\noi
\emph{Proof.} First we prove the consistency of the IMH by showing
that $M_d$ satisfies the IMH for some Turing degree $d$ in a forcing extension
of $V$.

Let $\kappa$ be Woodin with an inaccessible above in $V$. Let $G$ be
generic over $V$ for the L\'{e}vy 
collapse of $\kappa$ to $\omega$. Work now in $V[G]$. $\Sigma^1_2$
determinacy holds and, as there is still an inaccessible, $M_d$
exists for each Turing degree $d$.  
It follows that the theory of
$(M_d,\in)$ is constant on a cone of Turing degrees $d$. Let $d$ be a
Turing degree such that the theory of $(M_e,\in)$ is constant 
for Turing degrees $e$ at least that of $d$.

We claim that $M_d$, endowed with its definable classes, witnesses the
IMH. Indeed, suppose that $\varphi$ is a sentence true in some model
$M$ of height $\mu(d)$ compatible with $M_d$. By Jensen coding there
is a real $y$ such that $d$ is recursive in $y$, $\mu(y)=\mu(d)$ and $M$
is a definable inner model of $M_y$. Let $e$ be the Turing degree of
$y$. Then for some formula $\psi$, $M_e$ satisfies the sentence 

\bs

\noi
\emph{The inner model defined by $\psi$ (with some choice of parameters) 
satisfies $\varphi$.}

\bs

\noi
It follows that there is
an inner model of $M_d$
which satisfies $\varphi$, as desired. This proves the consistency of
the IMH.

\bs

To say that a countable $M$, together with its countable collection of
definable classes, satisfies IMH is simply a
$\Pi^1_1$-statement with a real coding $M$ as parameter, since one
only needs to quantify over outer models of $M$ of height $M\cap \Ord$.
Thus the assertion that there exists a Turing degree $d$ such that
$M_d$ (with its definable classes) satisfies 
IMH is a $\Sigma^1_2$-statement and hence absolute.
So the existence of a Woodin cardinal with an inaccessible
above implies that such an $M_d$ exists in $V$ (and indeed in $L$).

To prove the stronger statement that in $V$, $M_d$ satisfies the IMH for a
cone of $d$'s, one argues as follows. Say that a set of reals $X$ is
\emph{absolutely $\Delta^1_2$} iff there is a pair of $\Sigma^1_2$
formulas $\varphi(x)$, $\psi(x)$ such that $X$ consists of all
solutions to $\varphi(x)$ in $V$ and $\varphi$ is equivalent to the
negation of $\psi$ both in $V$ and all of its forcing extensions. 

\bs

\noi
Claim. Assume that there is a Woodin cardinal. Then determinacy holds
for absolutely $\Delta^1_2$ sets.

\bs

\noi
Proof of Claim. As before let $G$ be generic for the L\'{e}vy collapse
of the Woodin cardinal to $\omega$. Then $\Sigma^1_2$ determinacy
holds in $V[G]$. By the Moschovakis Third Periodicity theorem, 
(\cite{mosch} Theorem 6E.1), if $X$ is $\Sigma^1_2$ in $V[G]$ there is
a definable winning strategy in $V[G]$ for one of the players in the
game $G_X$. By the homogeneity of the L\'{e}vy collapse, it follows that
absolutely $\Delta^1_2$ sets are determined in $V$. This proves the
Claim.

\bs

As there is an inaccessible in $V$, $M_d$ exists for each Turing
degree $d$ in $V$. Now it follows from the Claim that in $V$,  for any
sentence $\varphi$, either for a cone of Turing degrees $d$,

$$M_d \models \varphi$$

\nod or for a cone of Turing degrees $d$,

$$M_d \models \neg\varphi,$$

\nod since the relevant games are absolutely $\Delta^1_2$.
Therefore in $V$ the theory of $(M_d,\in)$ is constant for
a cone of Turing degrees $d$. We can then apply the argument used 
earlier in $V[G]$ to conclude that also in $V$, 
$M_d$ satisfies IMH for $d$ in a cone of Turing degrees.
\hfill{$\Box$}\\

\centerline{\emph{Parameters and the Strong Inner Model Hypothesis}}

\bs

How can we introduce parameters into the Inner Model Hypothesis?
The following result shows that inconsistencies arise
without strong restrictions on the type of parameters allowed.

\begin{prop} (\cite{imh})
The Inner Model Hypothesis with arbitrary ordinal parameters or with arbitrary
real parameters is inconsistent.
\end{prop}

So instead we consider \emph{absolute parameters}, as in
\cite{stable}.
For any set $x$, the \emph{hereditary cardinality} of $x$, denoted
$\hcard(x)$, is the cardinality of the transitive closure of $x$.
If $V^*$ is an outer model of $V$, then a parameter $p$ is
\emph{absolute between $V$ and $V^*$} iff $V$ and $V^*$ have the same 
cardinals $\le\hcard(p)$ and some parameter-free formula
has $p$ as its unique solution in both $V$ and $V^*$.

\bs

\noi
\emph{Inner Model Hypothesis with locally absolute parameters} \ Suppose that
$p$ is absolute between $V$ and $V^*$
and $\varphi$ is a first-order sentence with
parameter $p$ which holds in
an inner model of $V^*$. Then $\varphi$ holds in an inner model of
$V$.

\bs

For a singular cardinal $\kappa$, a \emph{$\Box_\kappa$ sequence} is a 
sequence of the form $\lan C_\alpha\mid \alpha<\kappa^+$, $\alpha$ limit$\ran$
such that each $C_\alpha$ has ordertype less than $\kappa$ and for $\bar\alpha$
in $\Lim C_\alpha$, $C_{\bar\alpha}=C_\alpha\cap\bar\alpha$. \emph{Definable
$\Box_\kappa$} is the assertion that there exists a $\Box_\kappa$ sequence
which is definable over $H(\kappa^+)$ with parameter $\kappa$. We will
be interested in the special case $\kappa=\beth_\omega$, in which
case the parameter $\kappa$ is superfluous.

\begin{thm} \label{pimh}
The Inner Model Hypothesis with locally absolute parameters is
inconsistent. 
\end{thm}

\noi
\emph{Proof.} We first show that definable $\Box_\kappa$ fails, where
$\kappa$ is $\beth_\omega$. Let $\lan C_\alpha\mid 
\alpha<\kappa^+$, $\alpha$ limit$\ran$ be a $\Box_\kappa$ sequence
definable 
over $H(\kappa^+)$ without parameters. For each $n$ let $S_n$ be the
stationary set of all limit
$\alpha<\kappa^+$ such that the ordertype of $C_\alpha$ is greater
than  
$\beth_n$.

\bs

\noi
\emph{Claim.} Let $P_n$ be the forcing that adds a CUB subset of
$S_n$ using closed bounded subsets of $S_n$ as conditions, ordered by
end extension. Then $P_n$
is $\kappa^+$-distributive, i.e., does not add $\kappa$-sequences.

\bs

\noi
\emph{Proof of Claim.} It is enough to show that 
$P_n$ is $\beth_m^+$ distributive for each $m<\omega$. Assume that $m$
is  at least $n$. Suppose that $p$ is a condition and 
$\lan D_i\mid i<\beth_m\ran$ are dense. Let $\lan M_i\mid i<
\kappa^+\ran$ be a continuous chain of size $\kappa$ elementary
submodels  
of some large $H(\theta)$ such that $M_0$ contains $\kappa\cup
\{\lan C_\alpha\mid \alpha<\kappa^+$, $\alpha$ limit$\ran$, $p\}$ and
for each 
$i<\kappa^+$, $\lan M_j\mid j\le i\ran$ is an element of $M_{i+1}$.
Let $\kappa_i$ be $M_i\cap\kappa^+$ and $C$ the set of such
$\kappa_i$'s.  
Then a final segment $D$ of $C\cap \Lim C_\gamma$ is contained in
$S_n$, where  
$\gamma=\kappa_{\beth_m^+}$. Write $D$ as $\lan \kappa_{\alpha_i}\mid
i < \ordertype D\ran$. We can then choose a descending sequence 
$\lan p_i\mid i< \beth_m\ran$ of conditions below $p$ such that 
$p_{i+1}$ meets $D_i$ and belongs to $M_{\kappa_{\alpha_i}+1}$ for
each $i$.  
Then the greatest lower bound of this sequence meets each $D_i$. This 
proves the Claim.

\bs

It follows that for each $n$ the forcing $P_n$ does not alter 
$H(\kappa^+)$. By the Inner Model Hypothesis with locally absolute parameters
$S_n$ has a CUB subset $C_n$ in $V$ 
for each $n$. But this is a contradiction, as the intersection of the
$C_n$'s is empty.

\bs

Now we refine the above argument.
As not every real has a $\#$, there exist reals $R$ such that
$\kappa^+$ equals $\kappa^+$ of $L[R]$, where $\kappa$ is $\beth_\omega$.
Let $X$ be the set of such reals and for each $R$ in $X$ let 
$\lan C^R_\alpha\mid \alpha<\kappa^+$, $\alpha$ limit$\ran$ be the 
$L[R]$-least $\Box_\kappa$ sequence. Now for limit $\alpha<\kappa^+$,
define
$C^*_\alpha$ to be the intersection of the $C^R_\alpha$, $R\in X$. Then 
$\lan C^*_\alpha\mid\alpha<\kappa^+$, $\alpha$ limit$\ran$ is definable
in $H(\kappa^+)$ without parameters and has the 
properties of a $\Box_\kappa$ sequence with the sole exception that 
$C^*_\alpha$ is only guaranteed to be unbounded in $\alpha$ if $\alpha$ has
cofinality greater than $2^{\aleph_0}$. Now repeat the above argument
using $\lan C^*_\alpha\mid\alpha<\kappa^+$, $\alpha$ limit$\ran$ in place of 
$\lan C_\alpha\mid\alpha<\kappa^+$, $\alpha$ limit$\ran$, to obtain a
contradiction. $\Box$

\bs

To obtain the Strong Inner Model Hypothesis, we require
more absoluteness. We say that the parameter $p$ is 
\emph{(globally) absolute} iff there is a parameter-free formula which
has $p$ as 
its unique solution in all outer models of $V$ which have the same
cardinals $\le\hcard(p)$ as $V$.

\bs

\noi
\emph{Strong Inner Model Hypothesis (SIMH)} \ Suppose that $p$ is
absolute, $V^*$ is an outer model of $V$ with the same cardinals
$\le\hcard(p)$ as $V$ and $\varphi$ is a
first-order sentence with parameter $p$ 
which holds in an inner model of $V^*$. Then $\varphi$ holds in an
inner model of $V$.

\bs

\noi
\emph{Remark.} If above we assume that the sentence $\varphi$
holds not just in an inner model of $V^*$ but in $V^*$ itself, then in
the conclusion we may demand that in an inner model of $V$ witnessing
$\varphi$, $p$ is definable via the same
formula $\psi$ witnessing the absoluteness of $p$.
(This inner model may, however, fail to have the same cardinals
$\le\hcard(p)$ as $V$.) This is because we can replace the
sentence $\varphi$ by: 
``$\varphi$ holds and $p$ is defined by $\psi$''.

\begin{thm} \label{not-ch} (\cite{imh})
Assume the SIMH. Then CH is false. In fact, $2^{\aleph_0}$ cannot be
absolute and therefore cannot be 
$\aleph_\alpha$ for any ordinal $\alpha$ which is countable
in $L$.
\end{thm}

\begin{thm} \label{strong}
The SIMH implies the existence of an inner model with
a strong cardinal.
\end{thm}

\noi
\emph{Proof.} Assume not, and let $K$ be the core model below a strong
cardinal (see \cite{zeman}). As in the proof of Theorem \ref{meas}, we
let $K'$ 
denote the iterate of $K$ obtained by applying each order $0$ measure
exactly once. Then by the argument of Lemma \ref{k-prime-sing}, if
$\lambda$ is a 
cardinal 
then the $K'$-cofinality of $\lambda$ is not measurable in
$K'$. And by 
weak covering relative to $K$, if $\lambda$ is a singular cardinal,
then $\lambda^+$ is computed correctly in $K$ (i.e.,
$(\lambda^+)^K=\lambda^+$).

\begin{lem} \label{K'-correct} For any singular cardinal $\lambda$,
$\lambda^+$ is computed correctly in $K'$.
\end{lem}

\noi
\emph{Proof of Lemma \ref{K'-correct}.} 
This is clear if the $K$-cofinality of $\lambda$ is not measurable in
$K$, for then $\lambda$ 
is a fixed point of the iteration from $K$ to $K'$ and 
$(\lambda^+)^{K'}=(\lambda^+)^K =\lambda^+$. Otherwise let $\lan
K_i\mid i\in\Ord\ran$ result from the iteration of $K$
to $K'$ and choose $i$ so that the ultrapower map
$\sigma_i:K_i\to K_{i+1}$ applies the
order $0$ measure at $\kappa=\cof^{K_i}(\lambda)$. If $\lan 
\lambda_j\mid j<\kappa\ran$ is a continuous and increasing sequence in
$K_i$ with supremum $\lambda$, then $\lambda^+$ of $K_{i+1}$ is
represented in the ultrapower of $K_i$ by $\lan
\lambda_j^+\mid j<\kappa\ran$. In $K_i$, the product of the
$\lambda_j^+$'s contains a subset of size $(\lambda^+)^{K_i}$,
consisting of functions well-ordered by dominance on a final segment
of $\kappa$. It follows that $(\lambda^+)^{K_{i+1}}$ has cardinality
$\lambda^+$ and therefore $K_{i+1}$ computes $\lambda^+$ correctly. 
As $\lambda^+$ is a fixed point of the remaining iteration from
$K_{i+1}$ to $K'$, it follows that $K'$ computes $\lambda^+$
correctly. This proves Lemma \ref{K'-correct}.

\bs

We say that $\lambda$ is a \emph{cut point} of
$K'$ iff no extender on the $K'$ sequence with critical point less
than $\lambda$ has length at least $\lambda$ (i.e., $\lambda$ is
\emph{not overlapped} in $K'$). The class of cut points of $K'$ is
clearly closed. If the class of cut points of $K'$ is bounded, then
for sufficiently large $\lambda$ we can choose 
$f(\lambda)$ less than $\lambda$ which is the critical point of an
extender on the $K'$ sequence whose length is at least $\lambda$; but
then by Fodor's theorem, there would be a fixed $\kappa$ and extenders
on the $K'$ sequence with critical point $\kappa$ of arbitrarily large
length, which implies that $\kappa$ is a strong cardinal in $K'$. 
Thus as we have assumed that there is no
strong cardinal in $K'$, the class of cut
points of $K'$ is closed and unbounded.

\bs

Let $\lan\lambda_n\mid n\in\omega\ran$ be the first $\omega$-many
limit cardinals of $V$ which are cut points of $K'$, and let
$\lambda_\omega$ be their supremum. Then each $\lambda_n$, and of
course $\lambda_\omega$, has cofinality $\omega$.

\begin{lem} \label{absolute-parameter}
Each $\lambda_n$, and $\lambda_\omega$ as
well, is an absolute parameter.
\end{lem}

\noi
\emph{Proof of Lemma \ref{absolute-parameter}.} We first show that
$\lambda_0$ is absolute. Let $V^*$ be an 
outer model of $V$ with the same cardinals as $V$ up to
$\lambda_0$. Note that for some real $R$ in $V$, no $L_\alpha[R]$
satisfies ZFC and therefore $R^\#$ does not exist in $V^*$. It follows
that for any singular cardinal $\lambda$ of $V^*$, $\lambda$ is
singular in $V$ and $\lambda^+$ is computed correctly in $V$. In
particular, $V^*$ and $V$ have the same cardinals up to $\lambda_0^+$
and the same
singular cardinals up to $\lambda_0$.

\bs

It follows by Lemma \ref{K'-correct} that for any singular cardinal
$\lambda$ of $V^*$, $\lambda^+$ is computed 
correctly in both $K'$ and $(K^*)'$, where $(K^*)'$ denotes the 
$K'$ of $V^*$, obtained from $K^*$, the $K$ of $V^*$,
by applying each order $0$ measure exactly once. As both $K'$ and $(K^*)'$
are universal in $V^*$, it follows that they coiterate simply to
a common model $W$.

\bs

\noi
\emph{Claim.} 
The co-iteration of $K'$ with $(K^*)'$ fixes singular cardinals of
$V^*$ which are cut points either of $K'$ or of $(K^*)'$.

\bs

\noi
\emph{Proof of Claim.} 
Let $\lambda$ 
be a singular cardinal of
$V^*$ (and therefore also a singular cardinal of $V$). First assume that 
$\lambda$ is a cut point of $K'$. Then as $\lambda$ 
has non-measurable cofinality in 
$K'$, $\lambda$ is fixed by the iteration on the $K'$-side. And
as $\lambda$ has non-measurable cofinality in $(K^*)'$, 
$\lambda$ can only move on the $(K^*)'$-side if an extender overlapping 
$\lambda$ were applied. The assumption that $\lambda$ is not overlapped 
in $K'$ implies that $\lambda$ is not overlapped in $W$ 
and therefore the
least extender overlapping $\lambda$ was applied on the $(K^*)'$-side. 
But then $\lambda^+$ is not computed correctly in the resulting
ultrapower and therefore is computed correctly neither in $W$ nor in $K'$. 
This contradicts Lemma \ref{K'-correct}. 
As the same argument applies with $K'$ and $(K^*)'$ switched, this proves the 
Claim. 

\bs

It follows from the Claim that  $\lambda_0$ is the least 
limit cardinal of $V^*$ which is a cut point of $(K^*)'$. 
As $V^*$ is an arbitrary outer model of $V$ with the same
cardinals as $V$ up to $\lambda_0$, we have shown that $\lambda_0$ is
an absolute parameter.
The same argument shows that each $\lambda_n$
is absolute, and therefore so is $\lambda_\omega$, the supremum of the
first $\omega$ limit cardinals which are cut points of $K'$. This
proves Lemma \ref{absolute-parameter}.

\bs

Now let 
$\lan C_\alpha\mid\alpha<\lambda_\omega^+$, $\alpha$ limit$\ran$ be
the least 
$\Box_{\lambda_\omega}$ sequence of $K'$; this is also a $\Box_{\lambda_\omega}$ sequence
in $V$, as $(\lambda_\omega^+)^{K'}=\lambda_\omega^+$.
As in the proof of Theorem \ref{pimh}, 
there are generic extensions of $V$ preserving $H(\lambda_\omega^+)$ which add 
CUB subsets
to each $S_n=\{\alpha<\lambda_\omega^+\mid \ordertype
C_\alpha>\lambda_n\}$. It 
follows from the Strong Inner Model Hypothesis (and the Remark
immediately following its statement)
that for 
each $n$ there is an 
inner model $M_n$, with the correct $\lambda_\omega^+$ and
$\lambda_n$, in 
which $S_n^{M_n}$ 
contains a CUB subset $C_n$, where $S_n^{M_n}$ is defined 
using the least $\Box_{\lambda_\omega}$ sequence of $(K')^{M_n}$. The
latter may of course differ from the least $\Box_{\lambda_\omega}$
sequence of $K'$. However as $\lambda_\omega^+$ is computed correctly
in each $(K')^{M_n}$ and $\lambda_\omega$ is a cut point of
non-measurable cofinality in each $(K')^{M_n}$,
it follows that the 
$(K'\mid\lambda_\omega^+)^{M_n}$'s 
compare to a common $K''$ of height $\lambda_\omega^+$ with all 
ordinals in some CUB subset $C$ of $\lambda_\omega^+$ as closure
points. 
But if $\alpha$ is such a closure point in the intersection of the
$C_n$'s and $\alpha_n$ is the image of $\alpha$ 
under the comparison embedding of $(K'\mid\lambda_\omega^+)^{M_n}$ into $K''$,
then 
$C_{\alpha_n}$ as defined in $K''$ contains elements cofinal in
$\alpha$ and therefore $C_\alpha$ as defined in $K''$, an initial
segment of $C_{\alpha_n}$, has ordertype at least that of
$C_\alpha$ as defined in $(K'\mid \lambda_\omega^+)^{M_n}$. It follows
that 
$C_\alpha$ as defined in $K''$ has ordertype greater than $\lambda_n$
for each $n$, which is a contradiction. $\Box$

\bs

\noi
\emph{Remarks.} 
(a) It is likely that Theorem \ref{strong} can be improved to obtain
an inner model with a Woodin cardinal. But it is not possible to
obtain an iterable inner model with a Woodin cardinal and an
inaccessible above it (unless the SIMH is
inconsistent): Otherwise every real would be generic for Woodin's
extender algebra defined in an iterate of such an inner model,
implying that 
for every real $R$ there is an inaccessible in $L[R]$; this
contradicts Theorem \ref{imh}.\\
(b) David Asper\'o and the first author observed that
the consistency of the SIMH for the parameter $\omega_1$ follows as
in the proof of Theorem \ref{woodin2} from that of a Woodin cardinal
with an 
inaccessible above. In particular this yields the consistency of the
natural extension of L\'evy 
absoluteness asserting $\Sigma_1$ absoluteness with
parameter $\omega_1$ for arbitrary $\omega_1$-preserving 
extensions.

\bs

\noi
\emph{Question.} Is the Strong Inner Model Hypothesis
consistent relative to large cardinals?

\end{document}